\documentclass[11pt]{amsart}

\usepackage{amsmath}
\usepackage{amssymb}
\usepackage{amscd}
\usepackage{amsthm}
\usepackage{colortbl}
\usepackage{hyperref}
\usepackage{mathtools}
\usepackage{extarrows}

\usepackage[all]{xy}

\usepackage{tikz}
\usepackage{pgfplots}

\newtheorem{thm}{Theorem}[section]
\newtheorem{lemma}[thm]{Lemma}
\newtheorem{cor}[thm]{Corollary}

\newtheorem{prop}[thm]{Proposition}

\theoremstyle{remark}

\theoremstyle{definition}
\newtheorem{defn}[thm]{Definition}
\newtheorem{notation}[thm]{Notation}

\numberwithin{equation}{section}

\allowdisplaybreaks[4]

\begin{document}

\vfuzz0.5pc
\hfuzz0.5pc 

\newcommand{\claimref}[1]{Claim \ref{#1}}
\newcommand{\thmref}[1]{Theorem \ref{#1}}
\newcommand{\propref}[1]{Proposition \ref{#1}}
\newcommand{\lemref}[1]{Lemma \ref{#1}}
\newcommand{\coref}[1]{Corollary \ref{#1}}
\newcommand{\remref}[1]{Remark \ref{#1}}
\newcommand{\conjref}[1]{Conjecture \ref{#1}}
\newcommand{\questionref}[1]{Question \ref{#1}}
\newcommand{\defnref}[1]{Definition \ref{#1}}
\newcommand{\secref}[1]{\S \ref{#1}}
\newcommand{\ssecref}[1]{\ref{#1}}
\newcommand{\sssecref}[1]{\ref{#1}}

\newcommand{\RED}{{\mathrm{red}}}
\newcommand{\tors}{{\mathrm{tors}}}
\newcommand{\eq}{\Leftrightarrow}

\newcommand{\mapright}[1]{\smash{\mathop{\longrightarrow}\limits^{#1}}}
\newcommand{\mapleft}[1]{\smash{\mathop{\longleftarrow}\limits^{#1}}}
\newcommand{\mapdown}[1]{\Big\downarrow\rlap{$\vcenter{\hbox{$\scriptstyle#1$}}$}}
\newcommand{\smapdown}[1]{\downarrow\rlap{$\vcenter{\hbox{$\scriptstyle#1$}}$}}

\newcommand{\A}{{\mathbb A}}
\newcommand{\I}{{\mathcal I}}
\newcommand{\J}{{\mathcal J}}
\newcommand{\CO}{{\mathcal O}}

\newcommand{\C}{{\mathcal C}}
\newcommand{\BC}{{\mathbb C}}
\newcommand{\BQ}{{\mathbb Q}}
\newcommand{\m}{{\mathcal M}}
\newcommand{\h}{{\mathcal H}}
\newcommand{\Z}{{\mathcal Z}}
\newcommand{\BZ}{{\mathbb Z}}
\newcommand{\W}{{\mathcal W}}
\newcommand{\Y}{{\mathcal Y}}

\newcommand{\T}{{\mathcal T}}
\newcommand{\BP}{{\mathbb P}}
\newcommand{\CP}{{\mathcal P}}
\newcommand{\G}{{\mathbb G}}
\newcommand{\BR}{{\mathbb R}}
\newcommand{\D}{{\mathcal D}}
\newcommand{\DD}{{\mathcal D}}
\newcommand{\LL}{{\mathcal L}}
\newcommand{\f}{{\mathcal F}}
\newcommand{\E}{{\mathcal E}}
\newcommand{\BN}{{\mathbb N}}
\newcommand{\N}{{\mathcal N}}
\newcommand{\K}{{\mathcal K}}
\newcommand{\R} {{\mathbb R}}
\newcommand{\PP}{{\mathbb P}}
\newcommand{\Pp}{{\mathbb P}}
\newcommand{\BF}{{\mathbb F}}
\newcommand{\QQ}{{\mathcal Q}}
\newcommand{\closure}[1]{\overline{#1}}
\newcommand{\EQ}{\Leftrightarrow}
\newcommand{\imply}{\Rightarrow}
\newcommand{\isom}{\cong}
\newcommand{\embed}{\hookrightarrow}
\newcommand{\tensor}{\mathop{\otimes}}
\newcommand{\wt}[1]{{\widetilde{#1}}}
\newcommand{\ol}{\overline}
\newcommand{\ul}{\underline}

\newcommand{\bs}{{\backslash}}
\newcommand{\CS}{{\mathcal S}}
\newcommand{\CA}{{\mathcal A}}
\newcommand{\Q}{{\mathbb Q}}
\newcommand{\F}{{\mathcal F}}
\newcommand{\sing}{{\text{sing}}}
\newcommand{\U} {{\mathcal U}}
\newcommand{\B}{{\mathcal B}}
\newcommand{\X}{{\mathcal X}}

\newcommand{\ECS}[1]{E_{#1}(X)}
\newcommand{\CV}[2]{{\mathcal C}_{#1,#2}(X)}

\newcommand{\rank}{\mathop{\mathrm{rank}}\nolimits}
\newcommand{\codim}{\mathop{\mathrm{codim}}\nolimits}
\newcommand{\Ord}{\mathop{\mathrm{Ord}}\nolimits}
\newcommand{\Var}{\mathop{\mathrm{Var}}\nolimits}
\newcommand{\Ext}{\mathop{\mathrm{Ext}}\nolimits}
\newcommand{\EXT}{\mathop{{\mathcal E}\mathrm{xt}}\nolimits}
\newcommand{\Pic}{\mathop{\mathrm{Pic}}\nolimits}
\newcommand{\Spec}{\mathop{\mathrm{Spec}}\nolimits}
\newcommand{\Jac}{\mathop{\mathrm{Jac}}\nolimits}
\newcommand{\Div}{\mathop{\mathrm{Div}}\nolimits}
\newcommand{\sgn}{\mathop{\mathrm{sgn}}\nolimits}
\newcommand{\supp}{\mathop{\mathrm{supp}}\nolimits}
\newcommand{\Hom}{\mathop{\mathrm{Hom}}\nolimits}
\newcommand{\Sym}{\mathop{\mathrm{Sym}}\nolimits}
\newcommand{\nilrad}{\mathop{\mathrm{nilrad}}\nolimits}
\newcommand{\Ann}{\mathop{\mathrm{Ann}}\nolimits}
\newcommand{\Proj}{\mathop{\mathrm{Proj}}\nolimits}
\newcommand{\mult}{\mathop{\mathrm{mult}}\nolimits}
\newcommand{\Bs}{\mathop{\mathrm{Bs}}\nolimits}
\newcommand{\Span}{\mathop{\mathrm{Span}}\nolimits}
\newcommand{\IM}{\mathop{\mathrm{Im}}\nolimits}
\newcommand{\Hol}{\mathop{\mathrm{Hol}}\nolimits}
\newcommand{\End}{\mathop{\mathrm{End}}\nolimits}
\newcommand{\CH}{\mathop{\mathrm{CH}}\nolimits}
\newcommand{\Exec}{\mathop{\mathrm{Exec}}\nolimits}
\newcommand{\SPAN}{\mathop{\mathrm{span}}\nolimits}
\newcommand{\birat}{\mathop{\mathrm{birat}}\nolimits}
\newcommand{\cl}{\mathop{\mathrm{cl}}\nolimits}
\newcommand{\rat}{\mathop{\mathrm{rat}}\nolimits}
\newcommand{\Bir}{\mathop{\mathrm{Bir}}\nolimits}
\newcommand{\Rat}{\mathop{\mathrm{Rat}}\nolimits}
\newcommand{\aut}{\mathop{\mathrm{aut}}\nolimits}
\newcommand{\Aut}{\mathop{\mathrm{Aut}}\nolimits}
\newcommand{\eff}{\mathop{\mathrm{eff}}\nolimits}
\newcommand{\nef}{\mathop{\mathrm{nef}}\nolimits}
\newcommand{\amp}{\mathop{\mathrm{amp}}\nolimits}
\newcommand{\DIV}{\mathop{\mathrm{Div}}\nolimits}
\newcommand{\Bl}{\mathop{\mathrm{Bl}}\nolimits}
\newcommand{\Cox}{\mathop{\mathrm{Cox}}\nolimits}
\newcommand{\NE}{\mathop{\mathrm{NE}}\nolimits}
\newcommand{\NM}{\mathop{\mathrm{NM}}\nolimits}
\newcommand{\Gal}{\mathop{\mathrm{Gal}}\nolimits}
\newcommand{\coker}{\mathop{\mathrm{coker}}\nolimits}
\newcommand{\ch}{\mathop{\mathrm{ch}}\nolimits}

\def\O{\mathcal{O}}
\def\C{\mathbb{C}}
\def\CO{\mathcal{O}}
\def\P{\mathbb{P}}

\newcommand{\im}{\mathop{\mathrm{Im}}\nolimits}

\title{$\A^1$-curves on affine complete intersections}

\author{Xi Chen}
\address[Chen]{632 Central Academic Building\\
University of Alberta\\
Edmonton, Alberta T6G 2G1, CANADA}
\email{xichen@math.ualberta.ca}

\author{Yi Zhu}
\address[Zhu]{Pure Mathematics\\Univeristy of Waterloo\\Waterloo, ON N2L3G1\\ Canada}
\email{yi.zhu@uwaterloo.ca}


\date{\today}

\begin{abstract}
We generalize the results of Clemens, Ein, and Voisin regarding rational curves and zero cycles on generic projective complete intersections to the logarithmic setup. 
\end{abstract}

\maketitle
\section{Introduction}
In this paper, we work with varieties over the complex numbers. First we introduce the notion of \emph{smooth complete intersection pairs}.
\begin{defn}\label{not:ci}
	Let $X$ be a complete intersection in $\P^n$ of type $(d_1,\cdots,d_c)$. Let $D \subset X$ be a hypersurface section of degree $k$. We call the pair $(X,D)$ a {\em smooth complete intersection pair of type $(d_1,\cdots,d_c;k)$} if both $X$ and $D$ are smooth. We define the \emph{total degree} $d$ of the pair $(X,D)$ by $$d=d_1+\cdots+d_c+k.$$ When $k=0$, the boundary is empty and we simply denote $(X,D)$ by $X$.
	
\end{defn}


The existence of rational curves, algebraic hyperbolicity and rational equivalence of zero cycles on generic complete intersection of general type has been studied by Clemens \cite{Clemens86}, Ein \cite{Ein88,Ein91}, and Voisin \cite{Voisin94, Voisin96,Voisin98}.

\begin{thm}[Clemens,Ein,Voisin]\label{thm:CEV}
	Let $X$ be a generic complete intersection in $\P^n$ of type $(d_1,\cdots,d_c)$.
	\begin{enumerate}
		\item If $d\ge 2n-c$, $X$ has no rational curves;
		\item If $d\ge 2n-c+1$, $X$ is algebraically hyperbolic;
		\item If $d\ge 2n-c+2$, no two points of $X$ are rationally equivalent.
	\end{enumerate}
\end{thm}

The bounds above are not optimal. Voisin \cite{Voisin98} further improved the bound (1) to $d\ge 2n-2$ in case of hypersurfaces,  which is optimal because hypersurfaces of degree $\le 2n-3$ always admit lines. 

In this paper, we generalize Theorem \ref{thm:CEV} to smooth complete intersection pairs, where we study $\A^1$-curves and $\A^1$-equivalence of zero cycles instead. See Theorems \ref{thm:logEin}, \ref{thm:A1} and Corollary \ref{cor:chen} below. They specialize to Theorem \ref{thm:CEV} when the boundary is empty.

\subsection{$\A^1$-curves}
An \emph{$\A^1$-curve} is an algebraic map from $\A^1$ to the interior of a pair. When the boundary is empty, $\A^1$-curves are simply rational curves. We first study $\A^1$-curves on generic complete intersection pairs of general type. 
\begin{thm}\label{thm:logEin}
	Let $(X,D)$ be a generic complete intersection pair in $\P^n$ of type $(d_1,\cdots,d_c;k)$. If $d\ge 2n-c$, the interior $X-D$ contains no $\A^1$-curves.
\end{thm}

When the boundary is nonempty, the bound in Theorem \ref{thm:logEin} is optimal because a general such pair in $\P^n$ with $d\le 2n-c-1$ always admits an $\A^1$-line. Furthermore, we complete the last step studying $\A^1$-curves on complete intersection surface pairs in $\P^n$ of total degree $d$, summarized as the table below. 
\begin{center}
		\begin{tabular}{ | c | c | c | }
	\hline $\dim X=2$  & $(X,D)$ & 	$\A^1$-curves \\ \hline
	$d\le n$ & log Fano & log rationally connected \cite{CZ}  \\ \hline
	 $d=n+1$&  log K3& generically countable \cite{ChenX99,Li12,logK3} \\ \hline
	  $d\ge n+2$&of log general type& generically none (Thm. \ref{thm:logEin})\\
		\hline
		\end{tabular}
	\end{center}
	

\subsection{Algebraic hyperbolicity}
\begin{thm}\label{thm:logV}
	Let $(X,D)$ be a generic complete intersection pair in $\P^n$ of type $(d_1,\cdots,d_c;k)$. If $$d\ge 2n-c-l+1,$$ any closed subvariety $Y$ of $X-D$ of dimension $l$ has an effective log canonical bundle on its desingularisation; and if the equality is strict, $Y$ has a big log canonical bundle on its desingularisation. 
\end{thm}

Theorem \ref{thm:logV} implies algebraic hyperbolicity of such pairs. 

\begin{cor}\label{cor:chen}
	Let $X$ be a generic complete intersection in $\P^n$ of type $(d_1,\cdots,d_c;k)$. If $d\ge 2n-c+1$, the interior $X-D$ is algebraically hyperbolic.
\end{cor}

For generic complete intersection pairs of type $(1;k)$, Theorems \ref{thm:logEin}, \ref{thm:logV} and Corollary \ref{cor:chen} are proved by  the first named author \cite{ChenX-log} and Pacienza-Rousseau \cite{PR-log}.

\subsection{$\A^1$-equivalence of zero cycles}

For open varieties, the right substitution for Chow group of zero cycles is Suslin's $0$-th homology group $h_0(U)$, that is, the group of zero cycles modulo $\A^1$-equivalences. When the boundary is empty, it coincides with the Chow group of zero cycles. For surface pairs, the log version of Mumford's theorem and Bloch's conjecture was studied in \cite{logmumford,logbloch}. For generic complete intersection pairs, we have the following stronger version of Theorem \ref{thm:logEin}.

\begin{thm}\label{thm:A1}
	Let $(X,D)$ be a generic complete intersection pair in $\P^n$ of type $(d_1,\cdots,d_c;k)$. If $d\ge 2n-c+2$, no two points of the interior $X-D$ are $\A^1$-equivalent.
\end{thm}



\section{Global positivity}

In this section, we generalize Voisin's global positivity result \cite[Prop. 1.1]{Voisin96} for smooth complete intersection pairs. For the rest of the paper, we fix the following notations.

\begin{notation}
With the same notations as in Definition \ref{not:ci}, let $k:=d_{c+1}$. Let $S^{d_i}:=H^0(\O_{\P^n}(d_i))$ for $i=1,\cdots, c+1$. Let $S$ be the product vector space $\prod_i^{c+1} S^{d_i}$. We denote by $S^\circ$ the open subset of $S$ parametrizing smooth complete intersection pairs. Let $(\X,\D)\subset \P^n\times S^\circ$ be the universal family of smooth complete intersection pair. Let $\O_\X(1)$ be the pullback line bundle $pr^*_1(\O_{\P^n}(1))$. For any $t\in S^\circ$, denote by $(X_t,D_t)$ the smooth complete intersection pair parametrized by $t$. We assume that $\dim X_t\ge 2$. For any log pair $(Y,E)$, denote by $T_{Y^\dagger}$ its log tangent bundle $T_Y(-\log E)$.
\end{notation}


\begin{lemma}\label{lem:1}
	For any smooth complete intersection pair $(X,D)$ with $\dim X\ge 2$, we have $$H^0(\Omega^{i}_X(\log D))=0$$ for $0<i<\dim X$.
\end{lemma}
\proof The long exact sequence of the residue sequence gives 
$$ H^0(\Omega^{i}_X)\to H^0(\Omega^{i}_{X}(\log D))\to H^0(\Omega^{i-1}_{D})\to H^1(\Omega^{i+1}_X).$$ The first term vanishes by the Lefschetz hyperplane theorem. If $\dim D\ge 2$, the third term vanishes by the Lefschetz hyperplane theorem as well. If $\dim D=1$, the last map is the Gysin map which is injective. Therefore, $H^0(\Omega^{i}_{X}(\log D))=0$.\qed

\begin{lemma}\label{lem:2}
	If $d\ge n+2$, then  $h^0(T_{X_t^\dagger}(1))=0$ for every $t\in S^\circ$.
\end{lemma}

\proof By Serre duality and Lemma \ref{lem:1}, we have 
\begin{align*}
h^0(T_{X_t^\dagger}(1))&=h^0(\Omega^{n-1}_{X_t}(\log D_t)\otimes (-K_{X_t}-D_t)\otimes \O(1))\\
&=h^0(\Omega^{n-1}_{X_t}(\log D_t)\otimes \O(n+2-d))\\
&\le h^0(\Omega^{n-1}_{X_t}(\log D))=0.\qed
\end{align*}

\begin{prop}\label{prop:pos}
The log tangent bundle $$T\X^\dagger(1)|_{X_t}$$ is globally generated for every $t\in S^\circ $ if $h^0(T_{X_t^\dagger}(1))=0$.
\end{prop}

\proof 	By \cite[Lem. 4.1]{CZ}, we have the short exact sequence
\begin{equation*}
\xymatrix{	0 \ar[r] & \O_{\D}  \ar[r] & T\X^\dagger|_{\D}\ar[r]   & T_{\D} \ar[r]  & 0.	}
\end{equation*}
The global generation of $T\X^\dagger(1)|_{X_t}$ implies the global generation of $T\D(1)|_{D_t}$. In particular, Proposition \ref{prop:pos} for the nonempty boundary case implies the empty boundary case. For the rest of the proof, we assume that the boundary is nonempty.

Since $(\X,\D)$ is a log smooth family over $S^\circ$, we have
\begin{equation*}
\xymatrix{	0 \ar[r] & T_{X_t^\dagger}(1)  \ar[r] & T\X^\dagger(1)|_{X_t}\ar[r]   & S\otimes \O_{X_t}(1) \ar[r]  & 0.	}
\end{equation*}
By \cite[Lemma 2.1]{CZ}, the log tangent bundle $T\X^\dagger$ is determined by the short exact sequence:
\begin{equation*}
\xymatrix{	0 \ar[r] & T\X^\dagger \ar[r] & \O_{\X}(1)^{\oplus (n+1)}\ar[r]^-\alpha & \sum_{i=1}^{c+1} \O_\X(d_i)\ar[r] & 0,}
\end{equation*}
where $\alpha$ is given by the multiplication of the Jacobian. The above two sequences lead to the commutative diagram:

\begin{equation*}
\xymatrix@C-=0.5cm{	0 \ar[r] & T_{X_t^\dagger}(1) \ar[d]^{id} \ar[r] & T\X^\dagger(1)|_{X_t}\ar[r] \ar[d]  & S\otimes \O_{X_t}(1) \ar[r]\ar[d]^{ev}  & 0\\
0 \ar[r] & T_{X_t^\dagger}(1) \ar[r] & \O_{X_t}(2)^{\oplus (n+1)}\ar[r]^-\alpha & \sum_{i=1}^{c+1} \O_\X(d_i+1)|_{X_t}\ar[r] & 0.	}
\end{equation*}

Since $h^0(T_{X_t^\dagger}(1))=0$, we obtain the corresponding long exact sequences

\begin{equation*}
\xymatrix@C-=0.5cm{	0 \ar[r] & H^0(T\X^\dagger(1)|_{X_t}) \ar[d]^{} \ar[r] & S\otimes S^1\ar[rd]^\mu\ar[r]^\mu \ar[d]^{ev}  & H^1(T_{X_t^\dagger}(1)) \ar[r]^{}\ar[d]^{id}  & H^1(T\X^\dagger(1)|_{X_t}) \ar[d]\\
	0 \ar[r] & H^0(\O_{X_t}(2)^{\oplus (n+1)})\ar[r]^-\alpha & \prod_{i=1}^{c+1} S^{d_i+1}\ar[r] & H^1(T_{X_t^\dagger}(1) )\ar[r]^-\beta & H^1(\O_{X_t}(2)^{\oplus (n+1)}).	}
\end{equation*}

We have the following properties:
\begin{enumerate}
	\item $H^0(T\X^\dagger(1)|_{X_t}) =\ker(\mu)$;
	\item $\ker(\beta)= \prod_{i}^{c+1} S^{d_i+1}/\im(\alpha)$;
	\item since $ev$ is surjective, $\im(\mu)= \ker(\beta)$. Thus we have the map
	$$\mu: S\otimes S^1\to \ker(\beta).$$
\end{enumerate}

Now for any point $x\in X_t$, tensoring all the terms in the diagrams as above with the ideal sheaf $\I_x$, we have another commutative diagram:
\begin{equation*}
\xymatrix@C-=0.5cm{	 H^0(T\X^\dagger(1)|_{X_t}\otimes\I_x) \ar[d]^{} \ar[r] & S\otimes S^1_x\ar[rd]^{\mu_x}\ar[r]^-{\mu_x} \ar[d]^{ev_x}  & H^1(T_{X_t^\dagger}(1)\otimes\I_x) \ar[r]^{}\ar[d]^{id}  & H^1(T\X^\dagger(1)|_{X_t}\otimes\I_x) \ar[d]\\
  H^0(\O_{X_t}(2)^{\oplus (n+1)}\otimes\I_x)\ar[r]^-{\alpha_x} & \prod_{i}^{c+1} S^{d_i+1}_x\ar[r] & H^1(T_{X_t^\dagger}(1)\otimes\I_x )\ar[r]^-{\beta_x} & H^1(\O_{X_t}(2)^{\oplus (n+1)}\otimes\I_x),	}
\end{equation*}
where $S^m_x=H^0(\O_{X_t}(m)\otimes\I_x)$. We have the following properties:
\begin{enumerate}
	\item $H^0(T\X^\dagger(1)|_{X_t}\otimes\I_x) =\ker(\mu_x)$;
	\item $\ker(\beta_x)= \prod_{i}^{c+1} S^{d_i+1}_x/\im(\alpha_x)$;
	\item since $ev_x$ is surjective, we have $\im(\mu_x)= \ker(\beta_x)$. Thus we write 
	$$\mu_x: S\otimes S^1_x\to \ker(\beta_x).$$
\end{enumerate}

Finally, consider the commutative diagram:
$$
\xymatrix@C-=0.5cm{	0 \ar[r] \ar[d] & H^0(T_{X_t^\dagger}(1)|_{x}) \ar[rd]^\gamma\ar[d]^{} \ar[r] & H^1(T_{X_t^\dagger}(1)\otimes \I_x)\ar[r] \ar[d]^{\beta_x}  & H^1(T_{X_t^\dagger}(1)) \ar[r]^{}\ar[d]^{\beta}  & 0\\
	H^0(\O_{X_t}(2)^{\oplus (n+1)}) \ar[r]^{u} & H^0(\O_{X_t}(2)^{\oplus (n+1)}|_{x})\ar[r] &H^1(\O_{X_t}(2)^{\oplus (n+1)}\otimes \I_x)\ar[r] & H^1(\O_{X_t}(2)^{\oplus (n+1)})\ar[r]& 0.	}
$$

Since $\O_{X_t}(2)^{\oplus n+1}$ is globally generated, the map $u$ is surjective. In particular, the composite map $\gamma$ is the zero map. Hence we get 
\begin{equation*}
\xymatrix{	
	0 \ar[r] & H^0(T_{X_t^\dagger}(1)|_{x})  \ar[r] & \ker(\beta_x)\ar[r]   & \ker(\beta) \ar[r]  & 0.	}
\end{equation*}
It follows that $$\dim \im(\mu_x)-\dim\im(\mu)=\dim X_t.$$ On the other hand, we have 
$$\dim\ker(\mu)-\dim\ker(\mu_x)=\dim S\otimes S^1-\dim S\otimes S^1_x+\dim\im(\mu_x)-\dim\im(\mu)$$
$$=\dim S+\dim X_t=\dim \X.$$
Thus $$h^0(T\X^\dagger(1)|_{X_t})-h^0(T\X^\dagger(1)|_{X_t}\otimes\I_x)=\dim\X.$$ In particular, $T\X^\dagger(1)|_{X_t}$ is globally generated.\qed

Now Proposition \ref{prop:pos} implies 
\begin{cor}\label{cor:c}
	For any $l\ge 0$, the bundle $\wedge^l T\X^\dagger\otimes\O_{X_t}(l)$ is globally generated and the bundle $\wedge^l T\X^\dagger\otimes\O_{X_t}(l+1)$ is very ample if $d\ge n+2$.\qed
\end{cor}

\begin{cor}\label{cor:main}
	If $d\ge n+2$, then $\Omega_{\X^\dagger}^{\dim\X-l}|_{X_t}$ is globally generated when $d\ge l+n+1$ and is very ample when the inequality is strict. 
\end{cor}

\proof By Serre duality, we have 
\begin{align*}
\wedge^l T\X^\dagger\otimes\O_{X_t}(l)&= \Omega^{\dim\X-l}_{\X^\dagger}\otimes K_{\X^\dagger}^{-1}\otimes \O_{X_t}(l)\\&=\Omega^{\dim\X-l}_{\X^\dagger}\otimes \O_{X_t}(l+n+1-d).
\end{align*} Now the assertions follow from Corollary \ref{cor:c}.\qed

\section{Proof of Main Theorems}

\subsection{Proof of Theorems \ref{thm:logEin}, \ref{thm:logV}}
\proof[Proof of Theorem \ref{thm:logEin}] If not, there exists a family of $\A^1$-curves
$$
\xymatrix{
	V\times \A^1\ar[d]\ar[r]^f & (\X,\D)\ar[d]\\
	V\ar[r]^j &S^\circ,
	}
$$
where $j$ is an \'etale dominant morphism. By \cite[Lem. 3.1]{logmumford}, the morphism $f$ extends to a morphism of log pairs $$f:(V'\times \P^1, V'\times \{\infty\})\to (\X,\D),$$ where $V'$ is a dense open subset of $V$.
Here the dimension of $(V'\times \P^1, V'\times \{\infty\})$ is $\dim \X-(n-c)+1$. Corollary \ref{cor:main} with $l=n-c-1$ implies that	$\Omega_{\X^\dagger}^{\dim \X-l}|_{X_t}$ is globally generated if $$d\ge n-c-1+n+1=2n-c.$$ Pullback via $f$ gives a nontrivial section of the log canonical bundle of $(V'\times \P^1, V'\times \{\infty\})$, which is absurd because the pair is log uniruled. \qed

\proof[Proof of Theorem \ref{thm:logV}] The bound $d\ge 2n-c-l+1$ implies that $\Omega_{\X^\dagger}^{\dim \X-l}|_{X_t}$ is very ample. Now Theorem \ref{thm:logV} follows from the same proof as in \cite[Cor. 3]{PR-log}.\qed 

\subsection{$\A^1$-equivalence of two points}
To prove Theorem \ref{thm:A1}, we follow Voisin's approach \cite{Voisin94}. We first introduce a Mumford type invariant $\delta Z$ (Definition \ref{def}).
\begin{notation}
	Let $$\pi:(X,D)\to B$$ be a log smooth family of log pairs of relative dimension $n\ge 2$ with $\dim B=N$ and $D\neq \emptyset$. Assume that for every geometric fiber $X_b$,  $H^0(\Omega^i_{X_b}(\log D))=0$ for $0<i<n$. In particular, $H^1(\O_{X_b})=0$. Assume that there are two distinct sections $$p,q:B\to X-D$$ and denote the relative zero cycle $Z=p(B)-q(B)$. 
\end{notation}

The cycle $Z$ defines an $(n,n)$-form on $X$:$$[Z]\in H^n(\Omega^n_X).$$
Note that such form is supported on an open neighborhood of $Z$, in particular, it is away from $D$. We may assume that $$[Z]\in H^n(\Omega^n_X(\log D)\otimes \O_X(-D)).$$ Let $[Z_0]$ be the image of $[Z]$ via the natural map $$H^n(\Omega^n_X(\log D)\otimes \O_X(-D))\to H^0(B,R^n\pi_*(\Omega^n_X(\log D)\otimes \O_X(-D)).$$

Now by log smoothness, consider the short exact sequences
$$\xymatrix{
	0 \ar[r] & K\ar[r] & \Omega^n_X(\log D)\ar[r] & \Omega^n_{X|B}(\log D)\ar[r] & 0.}$$

After tensoring by $\O_X(-D)$, we get 
$$\xymatrix{
	0 \ar[r] & K(-D)\ar[r] & \Omega^n_X(\log D)(-D)\ar[r] & \Omega^n_{X|B}(\log D)(-D)\ar[r] & 0.}$$
Applying $R^n\pi_*$, since  $R^{n-1}\pi_*(\Omega^n_{X|B}(\log D)(-D))|_b$ is isomorphic to $H^1(\O_{X_b})$, which is zero, we have
$$\xymatrix{
	0 \ar[r] & R^n\pi_*K(-D)\ar[r] & R^n\pi_*\Omega^n_X(\log D)(-D)\ar[r] & R^n\pi_*\Omega^n_{X|B}(\log D)(-D). }$$ On the other hand, since $[Z_0]$ is a homologically trivial zero cycle, $[Z_0]$ maps to zero in $R^n\pi_*\Omega^n_{X|B}(\log D)(-D)$. 

\begin{defn}\label{def}
We define $\delta Z:=[Z_0]\in H^0(R^n\pi_*(K(-D))).$
\end{defn}

We now give a concrete description of $\delta Z$. By Serre duality, we have
	\begin{align*}
	R^n\pi_*(\Omega^n_X(\log D)(-D))_b&\cong H^n(\Omega_{X}^n(\log D)(-D)|_{X_b})\\
	&\cong H^0((\Omega_{X}^n(\log D)(-D))^\vee|_{X_b}\otimes K_{X_b})^*\\
	&\cong H^0((\Omega_{X}^N(\log D)\otimes K_X)|_{X_b}\otimes K_{X_b})^*\\
	&\cong H^0(\Omega_{X}^N(\log D)\otimes \pi^*K_B^{-1}|_{X_b})^*,
	\end{align*}
	and 
	\begin{align*}
R^n\pi_*\Omega^n_{X|B}(\log D)(-D)|_b\cong H^n(\Omega^n_{X_b}(\log D)(-D)) \cong H^0(\O_{X_b})\cong \C.
	\end{align*} 
	
	Hence we have $$R^n\pi_*K(-D)|_b\cong H^n(K(-D)|_{X_b})\cong (H^0(\Omega_{X}^N(\log D)\otimes \pi^*K_B^{-1}|_{X_b})/H^0(\O_{X_b}))^*$$ and the class $\delta Z|_b$ gives an element in $$(H^0(\Omega_{X}^N(\log D)\otimes \pi^*K_B^{-1}|_{X_b})/H^0(\O_{X_b}))^*.$$

\begin{lemma}\label{lem:distinct}
	The invariant $\delta Z$ is the pullback $p^*-q^* $ given as below:
	$$H^0(\Omega_{X}^N(\log D)\otimes \pi^*K_B^{-1}|_{X_b})\to H^0(\Omega_{B}^N(\log D)\otimes K_B^{-1}|_{X_b})\cong \C$$
\end{lemma}

\proof Since both sections $p$ and $q$ are away from the boundary, the proof directly follows from \cite[Prop 1.16]{Voisin94}.\qed

Next consider the short exact sequence
$$\xymatrix{
	0 \ar[r] & \pi^*\Omega^n_B\ar[r] & K\ar[r] & Q\ar[r] & 0,}$$ where $Q$ admits a natural decreasing filtration whose graded pieces are $ \Omega^i_{X|B}(\log D)\otimes \pi^*\Omega^{n-i}_{B}$, for $0<i<n$. We have 
\begin{equation}\label{seq:dec}
\xymatrix{
	R^{n-1}\pi_*Q \ar[r]^-\psi & \Omega^n_B\otimes R^n\pi_*\O_X(-D)\ar[r] & R^n\pi_*(K(-D))\ar[r] & R^n\pi_*Q(-D).}
\end{equation}

\begin{lemma}\label{lem:van} If for any $0<i<n$, $H^0(\Omega^i_{X_b}(\log D))=0$, then 
	$$R^n\pi_*Q(-D)=0.$$
\end{lemma}

\proof It suffices to verify for the graded pieces of $Q$. For each $0<i<n$, we have \begin{align*}
R^n\pi_* \Omega^i_{X|B}(\log D)(-D)\otimes \Omega^{n-i}_{B}|_b&\cong H^n(\Omega^i_{X_b}(\log D)(-D)) \otimes \Omega^{n-i}_{B}|_b\\
&\cong H^0(\Omega^i_{X_b}(\log D)) \otimes \Omega^{n-i}_{B}|_b=0.
\end{align*}\qed

Lemma \ref{lem:van} and (\ref{seq:dec}) imply

\begin{lemma}\label{lem:3}
	$\delta Z\in H^0(\Omega^n_B\otimes R^n\pi_*\O_X(-D)/\im\psi)$.\qed
\end{lemma}

\begin{lemma}\label{lem:0}
	If $[Z]$ is $\A^1$-equivalent to $0$, then $\delta Z=0$ on some open subset of $B$.
\end{lemma}

\proof Lemma \ref{lem:3} implies that locally, $\delta Z$ is an $(n,n)$-form decomposed as $\Omega^n_B\otimes R^n\pi_*\O_X(-D)$ and vertically compactly supported on the tubular neighborhood $U$ of $Z$. Thus $\delta Z$ is equivalent to the morphism
\begin{align*}
\delta Z:H^0(\pi_*\Omega^n_{X|B}(\log D))&\to H^0(\Omega^n_B)\\
\omega&\mapsto  \int_{} \omega\wedge \delta Z
\end{align*} by taking the wedge product and integrating along the fiber of $\pi$.

For any $\omega\in H^0(\pi_*\Omega^n_{X|B}(\log D))$, we know that $$\int_Z \omega|_{Z}\wedge \eta=\int_X \omega\wedge \delta Z\wedge \eta $$ for any compactly supported closed form $\eta$. By the relative version of \cite[Theorem 3.3]{logmumford}, $\omega|_{Z}=0$. (\cite[Theorem 3.3]{logmumford} is stated for surface pairs, but the proof works for any pairs.) Thus the integration above is zero. We conclude that the form $\omega\wedge \delta Z$ is exact. It is also vertically compactly supported on the tubular neighborhood $U$.

Now the Thom isomorphism tells us integration along the vertical fiber of an exact form gives us an zero element in $H^n(B)$. Therefore, $\delta Z(\omega)=0$. \qed





\proof[Proof of Theorem \ref{thm:A1}] If the conclusion is not true, there exists a smooth variety $S'$ \'etale dominant over an open subset of $S$ such that 
\begin{itemize}
	\item the base change $(\X,\D)\times_S S'$ admits two distinct sections $p$ and $q$ over $S'$;
	\item the relative zero cycle $Z=p(S')-q(S')$ is trivial under $\A^1$-equivalence.
\end{itemize}

By Lemma \ref{lem:0}, $\delta Z$ vanishes on some open subset of $S'$. On the other hand, by Corollary \ref{cor:main}, $\Omega_{\X^\dagger}^{\dim\X-\dim X_t}|_{X_t}$ is very ample if $$d\ge n-c+n+1+1=2n-c+2.$$
This implies that $p^*-q^*$ is not zero over a general point. Thus by Lemma \ref{lem:distinct}, $\delta Z$ is not zero. We have a contradiction.
\qed


\bibliography{myref}	

\begin{thebibliography}{Zhu15}

\bibitem[Che99]{ChenX99}
Xi~Chen.
\newblock Rational curves on {$K3$} surfaces.
\newblock {\em J. Algebraic Geom.}, 8(2):245--278, 1999.

\bibitem[Che04]{ChenX-log}
Xi~Chen.
\newblock On algebraic hyperbolicity of log varieties.
\newblock {\em Commun. Contemp. Math.}, 6(4):513--559, 2004.

\bibitem[Cle86]{Clemens86}
Herbert Clemens.
\newblock Curves on generic hypersurfaces.
\newblock {\em Ann. Sci. \'Ecole Norm. Sup. (4)}, 19(4):629--636, 1986.

\bibitem[CZ14]{CZ}
Qile Chen and Yi~Zhu.
\newblock Very free curves on {F}ano complete intersections.
\newblock {\em Algebr. Geom.}, 1(5):558--572, 2014.

\bibitem[CZ16]{logK3}
Xi~Chen and Yi~Zhu.
\newblock $\mathbb{A}^1$-curves on log {K}3 surfaces.
\newblock 2016.
\newblock Submitted. arXiv:1602.02204.

\bibitem[Ein88]{Ein88}
Lawrence Ein.
\newblock Subvarieties of generic complete intersections.
\newblock {\em Invent. Math.}, 94(1):163--169, 1988.

\bibitem[Ein91]{Ein91}
Lawrence Ein.
\newblock Subvarieties of generic complete intersections. {II}.
\newblock {\em Math. Ann.}, 289(3):465--471, 1991.

\bibitem[LL12]{Li12}
Jun Li and Christian Liedtke.
\newblock Rational curves on {K}3 surfaces.
\newblock {\em Invent. Math.}, 188(3):713--727, 2012.

\bibitem[PR07]{PR-log}
Gianluca Pacienza and Erwan Rousseau.
\newblock On the logarithmic {K}obayashi conjecture.
\newblock {\em J. Reine Angew. Math.}, 611:221--235, 2007.

\bibitem[Voi94]{Voisin94}
Claire Voisin.
\newblock Variations de structure de {H}odge et z\'ero-cycles sur les surfaces
  g\'en\'erales.
\newblock {\em Math. Ann.}, 299(1):77--103, 1994.

\bibitem[Voi96]{Voisin96}
Claire Voisin.
\newblock On a conjecture of {C}lemens on rational curves on hypersurfaces.
\newblock {\em J. Differential Geom.}, 44(1):200--213, 1996.

\bibitem[Voi98]{Voisin98}
Claire Voisin.
\newblock A correction: ``{O}n a conjecture of {C}lemens on rational curves on
  hypersurfaces'' [{J}.\ {D}ifferential {G}eom.\ {\bf 44} (1996), no.\ 1,
  200--213; {MR}1420353 (97j:14047)].
\newblock {\em J. Differential Geom.}, 49(3):601--611, 1998.

\bibitem[YZ15]{logbloch}
Qizheng Yin and Yi~Zhu.
\newblock $\mathbb{A}^1$-equivalence of zero cycles on surfaces {II}.
\newblock 2015.
\newblock Submitted. arXiv:1512.09079.

\bibitem[Zhu15]{logmumford}
Yi~Zhu.
\newblock $\mathbb{A}^1$-equivalence of zero cycles on surfaces.
\newblock 2015.
\newblock Submitted. arXiv:1510.01712.

\end{thebibliography}
\bibliographystyle{alpha}

\end{document}